\documentclass[a4paper,11pt]{article}
\usepackage{palatino}
\usepackage{textcomp}
\usepackage{amssymb}
\usepackage{graphicx}
\usepackage{epsfig}
\usepackage{epstopdf}
\usepackage{tikz}
\usepackage{mathrsfs}
\usepackage{picins}
\usepackage[english]{babel}
\usepackage{amscd}
\usepackage{amsthm}
\usepackage{amsmath}
\usepackage{fancyhdr}

\newtheorem{thm}{Theorem}[section]
\newtheorem{lemma}[thm]{Lemma}

\newtheorem{cor}[thm]{Corollary}

\newtheorem{conj}[thm]{Conjecture}

\renewcommand{\leq}{\leqslant}
\renewcommand{\geq}{\geqslant}
\renewcommand{\setminus}{\smallsetminus}
\author{Christos Pelekis\thanks{Delft University, Institute of Applied Mathematics, Mekelweg 4, 2628CD, Delft, Netherlands,  \texttt{pelekis.chr@gmail.com}}}

\begin{document}

\title{Bernoulli trials of fixed parity, random and randomly oriented graphs}

\maketitle

\abstract{Suppose you can color $n$ \emph{biased} coins with $n$ colors, all coins having the same bias.
It is forbidden to color both sides of a coin with the same color, but all other colors are allowed.
Let $X$ be the number of different colors after a toss of the coins. We present a method to obtain
an upper bound on a median of $X$. Our method is based on
the analysis of the probability distribution of
the number of vertices with even in-degree in  graphs whose
edges are given random orientations. Our analysis applies to the
distribution of the number of vertices with odd degree in random sub-graphs of
fixed graphs.
It turns out that there are parity restrictions on the random variables that are under consideration.
Hence, in order to present
our result,
we introduce a class of Bernoulli random variables whose total number of successes is of fixed parity and
are closely related to Poisson trials
conditional on the event that their outcomes have fixed parity.}

\section{Introduction}

The main motivation behind this work is the following problem that arose in the analysis
of a network coloring game (see \cite{PelSch}).
Suppose you can color $n$ \emph{biased} coins with $n$ colors, all coins having the same bias.
It is forbidden to color both sides of a coin with the same color, but all other colors are allowed.
Let $X$ be the number of different colors after a toss of the coins.
In what way should you color the coins in order to maximize the median of $X$?
What about upper bounds on the median of $X$?
In this paper we focus on the second question.
In previous work  (see \cite{PelSch}) we presented a method to obtain upper bounds on the median of $X$ in the case
of fair coins. In this work we extend this method to the case of biased coins.
Our analysis is heavily based on  the following model. \\
Suppose that  $G=(V,E)$ is a \emph{connected} graph on $n$ vertices and $m\geq n-1$ edges.
For every edge $e\in E$ call arbitrarily one of its endpoints \emph{head} and the other endpoint \emph{tail},
and consider the model in which each edge is getting a random orientation, independently and
is either oriented towards its tail with probability $p$ or
oriented towards its head with probability $1-p$.
The case in which all edges are oriented equiprobably towards one of the two possible directions
has been well studied (see for example \cite{AlmLinusson},\cite{Diarmir},\cite{Linusson}). Let $E_G$ be the
number of vertices of even in-degree
after assigning a random orientation on the edges of $G$.
In  \cite{PelSch} we computed the
distribution of $E_G$ in the case where $p$ equals
$\frac{1}{2}$. In this paper we  study the distribution of $E_G$ in the general case where $p\in (0,1)$.
We present a method to estimate the probability distribution of $E_G$ from below, in the
sense of stochastic orderings.
It turns out that $E_G$ has the same parity as $m-n$, i.e.,
$E_G = m-n\; \text{mod}\; 2$, a fact that
imposes parity restrictions on the number of vertices with even in-degree. \\
Our method applies to the distribution of the number of vertices with \emph{odd} degree, $O_{n,p}(G)$,
in random sub-graphs of a fixed graph, $G$, on $n$ vertices in which we either erase an edge with probability $1-p$
or keep it with probability $p$, independently for all edges.
Again, the degree-sum formula imposes parity restrictions on $O_{n,p}(G)$.
In particular $O_{n,p}(G)$ has to be even.  \\
Thus, in order to present our results,  we begin
by defining a class of Bernoulli random variables whose total number of successes is
of fixed parity and does not seem to have been studied before.
This class contains Bernoulli random variables that
are closely related to Poisson trials conditional on the event that
their outcome is of fixed parity. We study this class  in
Section \ref{FixPar}. In
Section \ref{randorient} we prove that $E_G$ is stochastically larger than
a certain random variable from this class and then we use this result to obtain an upper bound
on the median of the number of different colors after a toss of colored coins.
In Section \ref{randgraphsection} we apply our method to obtain a result on the distribution of the number of
vertices with odd degree in random graphs.
In Section \ref{apll} we employ our results to obtain probabilistic proofs of known results from the literature.
Finally, in Section \ref{open}, we consider some open questions.

\section{Bernoulli Trials of Fixed Parity}\label{FixPar}

In this section we define and state basic properties of a class of discrete probability distributions that
arise in the analysis of the random variables that are under consideration.
We study Bernoulli random variables conditioned on the total number of successes
having fixed parity. There has
been quite some work on Bernoulli random variables conditioned on the total number of successes
being at least a certain given value (see \cite{Meester} and references therein).
We begin by fixing some definitions and notation. \\

Denote by $B(n,p)$ a binomially distributed
random variable of parameters $n$ and $p$. That is, $B(n,p)$ is the number of
successes in $n$ independent and identical Bernoulli trials, $\text{Ber}(p)$.
A random variable that generalizes the binomial
is defined in the following way.
Fix a set of $n$ parameters, $I=\{p_1,\ldots,p_n\}$, from $(0,1)$ we denote by $\mathcal{H}(I)$
the random variable that counts the number of successes in $n$ independent, non-identical Bernoulli trials,
$\text{Ber}(p_i),i=1\ldots,n$. In other
words,
$\mathcal{H}(I)$ counts the number of $1$'s after a toss of
$n$ independent $0/1$-coins, $c_i,i=1,\ldots,n$,
having the property that coin $c_i$ shows $1$, or is a \emph{success}, with probability $p_i$.
The distribution of $\mathcal{H}(I)$ is well studied and is referred to as
\emph{Poisson binomial distribution}, or as \emph{Poisson trials}
(see \cite{Hoeffding},\cite{Wang}).
Our first result, concerning the parity of such a random variable, will be used repeatedly.

\begin{lemma}
\label{binparity}
Let $I=\{p_1,\ldots,p_n\}$
and $h_n := \mathcal{H}(I) \; \text{mod}\;2$. Then $h_n$ is a biased $0/1$ coin that
shows $1$ with probability $\frac{1}{2}(1-\prod_{i=1}^{n}(1-2p_i))$. That is, the probability that a
$\mathcal{H}(I)$ random variable is even equals
\[\frac{1}{2}(1+\prod_{i=1}^{n}(1-2p_i)) \]
and the probability that it is odd equals
\[\frac{1}{2}(1-\prod_{i=1}^{n}(1-2p_i)) . \]
\end{lemma}
\begin{proof}
The proof is by induction on $n$. When $n=1$ the conclusion is true. Suppose that it is true
when $|I|=n-1$ and consider a $\mathcal{H}(I)$ random variable with $|I|=n$. Then
\begin{eqnarray*}
 \mathbb{P}[\mathcal{H}(I) \; \text{even}] &=& p_n \cdot\mathbb{P}[\mathcal{H}(I\setminus\{p_n\}) \; \text{odd}]\\
  &+& (1-p_n)\cdot \mathbb{P}[\mathcal{H}(I\setminus\{p_n\}) \; \text{even}] \\
 &=& p_n \cdot\frac{1}{2}(1-\prod_{i=1}^{n-1}(1-2p_i)) \\
 &+& (1-p_n)\cdot \frac{1}{2}(1+\prod_{i=1}^{n-1}(1-2p_i)) \\
 &=&  \frac{1}{2} (1+\prod_{i=1}^{n}(1-2p_i))
\end{eqnarray*}
Since $\mathbb{P}[\mathcal{H}(I) \; \text{odd}]=1-\mathbb{P}[\mathcal{H}(I) \; \text{even}]$ the lemma follows.
\end{proof}

Given a set of parameters $I=\{p_1,\ldots,p_n\}$, set
$\alpha(I):= \mathbb{P}[\mathcal{H}(I)\; \text{even}]$
and $\beta(I)=1-\alpha(I)$.\\

Now fix a set of parameters $I=\{p_1,\ldots,p_n\}$ and
define a random variable whose outcomes have fixed parity, in the following way.
First consider the case of even outcomes.
Place the $0/1$ coins $c_1,\ldots ,c_n$ on a line. Roll a \emph{biased} die with $n$ faces
that shows $i \in \{1,\ldots, n\}$ with probability $\pi_i$.
That is, let $\pi=(\pi_1,\ldots, \pi_n)$ be such that $\sum \pi_i =1$ and
choose $i$ with probability $\pi_i$.
If the result of the die is $i\in\{1,\ldots,n\}$,
then \emph{toss all coins except} $c_i$. If the outcome after the toss has an even number of $1$'s,
then fix the parity by letting $c_i$ to be $0$. If the outcome has an odd
number of $1$'s, then
fix the parity by letting $c_i$ to be $1$. The number of $1$'s that we see after this (slightly dependent) toss is random. Denote it by $\mathcal{E}(I,\pi)$ and call this dependent toss an \emph{even-sum} toss of $n$ coins.
Similarly we define the \emph{odd-sum} toss of $n$ coins and denote by $\mathcal{O}(I,\pi)$
the number of $1$'s that we see after an odd-sum toss of $n$ coins.
Formally, for an even $k$, the probability distribution $\mathcal{E}(I,\pi)$ is defined by
\[ \mathbb{P}[\mathcal{E}(I,\pi)=k]= \sum_{i=1}^{n} \pi_i \cdot\left\{
 \mathbb{P}[\mathcal{H}(I\setminus\{p_i\})=k]+\mathbb{P}[\mathcal{H}(I\setminus\{p_i\})=k-1]\right\}  \]
and similarly for an odd $\ell$, the distribution of $\mathcal{O}(I,\pi)$ is defined by
\[ \mathbb{P}[\mathcal{O}(I,\pi)=\ell]=  \sum_{i=1}^{n}\pi_i \cdot\left\{
 \mathbb{P}[\mathcal{H}(I\setminus\{p_i\})=\ell]+\mathbb{P}[\mathcal{H}(I\setminus\{p_i\})=\ell-1] \right\}  \]

Note that in case all parameters $p_i\in I,i=1,\ldots, n$, are equal to $p$, then the probability
distribution of an even-sum toss equals
\[ \mathbb{P}[\mathcal{E}(I,\pi)=k] = \mathbb{P}[B(n-1,p)=k] + \mathbb{P}[B(n-1,p)=k-1] , \]
and so does not dependent on the vector $\pi=(\pi_1,\ldots, \pi_n)$. Similarly for the odd-sum toss.
In case all parameters $p_i$ are equal to $p$ we will
denote the random variables that count the number of successes in an even-sum (resp. odd-sum) toss of $n$ coins
by $A(n,p)$ (resp. $P(n,p)$).   \\

Notice also that in case $p_i=\frac{1}{2}$, for all $i\in \{1,\ldots,n\}$, the above formulas reduce to
\begin{eqnarray*} \mathbb{P}[A(n,1/2)=k]&=& \mathbb{P}\left[B(n-1,1/2)=k\right]+\mathbb{P}\left[B(n-1,1/2)=k-1\right] \\
&=& \binom{n-1}{k} \frac{1}{2^{n-1}} +  \binom{n-1}{k-1} \frac{1}{2^{n-1}} \\
&=& \binom{n}{k} \frac{1}{2^{n-1}}
\end{eqnarray*}
and similarly for $P(n,1/2)$.

The random variables just defined are related to
the random variable $\mathcal{H}(I)$, conditional on the event
that its outcomes have fixed parity.
More precisely, denote by $\mathcal{H}(I,0)$
(resp.$\mathcal{H}(I,1)$) the  random variable that has the same
distribution as $\mathcal{H}(I)$ conditional on the event that it's total number of successes
is even (resp. odd). That is, for even $k$
\[ \mathbb{P}[\mathcal{H}(I,0)=k] = \frac{1}{\alpha(I)} \mathbb{P}[\mathcal{H}(I)=k] \]
and, for an odd $\ell$,
\[ \mathbb{P}[\mathcal{H}(I,1)=\ell] = \frac{1}{\beta(I)} \mathbb{P}[\mathcal{H}(I)=\ell] . \]

Hence we can obtain an outcome of a, say, $\mathcal{H}(I,0)$ random
variable by tossing the coins again and again until we see an even outcome.
In case $I$ consists of $n$ parameters all equal to $p$, we will write $B(n,p,0)$ for $\mathcal{H}(I,0)$ and
$B(n,p,1)$ for $\mathcal{H}(I,1)$. Thus $B(n,p,0)$ is the random variable whose distribution function
is binomial, conditional on the event that the outcomes are even. Similarly for $B(n,p,1)$.

The following results shows the relation between conditional Poison trials and
the Bernoulli random variables that are under consideration.

\begin{lemma}
\label{evendistr}
If $I=\{p_1,\ldots,p_n\}$ and $\pi=(\pi_1,\ldots, \pi_n)$ is a probability vector
then the distribution of $\mathcal{E}(I,\pi)$ is the same as the distribution of the random variable that
takes even outcomes according to the following procedure. Roll a biased die with $n$ faces. If the result of the die is $i\in\{1,\ldots,n\}$ with probability $\pi_i$, then toss a $0/1$ coin having probability of showing $1$ equal to $1-\alpha(I\setminus\{p_i\})=\beta(I\setminus\{p_i\})$. If the outcome of this coin is $0$, then draw from a $\mathcal{H}(I\setminus\{p_i\},0)$ random variable  and add $0$. If the outcome is $1$, then draw from a $\mathcal{H}(I\setminus \{p_i\},1)$ random variable and add $1$.
\end{lemma}
\begin{proof} For an even $k$, write
\begin{eqnarray*}
\mathbb{P}[\mathcal{E}(I,\pi)=k]&=& \sum_{i=1}^{n} \pi_i \cdot \left\{
\alpha(I\setminus\{p_i\})\cdot \frac{\mathbb{P}[\mathcal{H}(I\setminus\{p_i\})=k]}{\alpha(I\setminus\{p_i\})}\right\} \\
&+&  \sum_{i=1}^{n} \pi_i \cdot \left\{ \beta(I\setminus\{p_i\})\cdot\frac{\mathbb{P}[\mathcal{H}(I\setminus\{p_i\})=k-1]}{\beta(I\setminus\{p_i\})} \right\} ,
\end{eqnarray*}
which can be rewritten as
\begin{eqnarray*}
\mathbb{P}[\mathcal{E}(I,\pi)=k]&=& \sum_{i=1}^{n} \pi_i \cdot
 \alpha(I\setminus\{p_i\})\cdot \mathbb{P}[\mathcal{H}(I\setminus\{p_i\},0)=k]  \\
 &+& \sum_{i=1}^{n} \pi_i \cdot  \beta(I\setminus\{p_i\})\cdot\mathbb{P}[\mathcal{H}(I\setminus\{p_i\},1)=k-1]
\end{eqnarray*}
and finishes the proof of the lemma.
\end{proof}

For random variables $Y, W$ that take values on the same sets, we will write $Y\sim W$
whenever $Y$ and $W$ have the same distribution.
%Given $x\in(0,1)$, denote by $\gamma_x$ the outcome of a $0/1$ coin that shows up $1$ with probability $x$.
%That is $\gamma_x\in \{0,1\}$.
Note that, in case all parameters $p_i$ are equal to $p$, the previous lemma says that
$A(n,p)$ has the same distribution as the random variable that takes even outcomes according to the
following procedure. Toss a $0/1$ coin whose probability of showing $1$ equals $\beta(\{p\}_{n-1})$. If the outcome
is a $1$, then toss $n$ independent $0/1$ coins that show up $1$ with probability $1$ until you see
an odd outcome, and add a $1$. If the outcome is $0$, then toss $n$ independent $0/1$ coins that show $1$ with probability  $p$ until you see an even outcome, and add a $0$ to this outcome.
We can formally express this as
\[ A(n,p) \sim B(1,\beta(\{p\}_{n-1})) + B(n-1,p, B(1,\beta(\{p\}_{n-1}))) . \]

Similarly one can prove the following result for $\mathcal{O}(I,\pi)$.

\begin{lemma}
\label{odddistr}
If $I=\{p_1,\ldots,p_n\}$ and $\pi=(\pi_1,\ldots, \pi_n)$ is a probability vector,
then the distribution of $\mathcal{O}(I,\pi)$ is the same as the distribution of the random variable that
takes odd outcomes according to the following procedure. Roll a biased die with $n$ faces. If the result of the die is $i\in\{1,\ldots,n\}$ with probability $\pi_i$, then toss a $0/1$ coin having probability of showing $1$ equal to $\alpha(I\setminus\{p_i\})=1-\beta(I\setminus\{p_i\})$. If the outcome of the coin is a $0$, then draw from a $\mathcal{H}(I\setminus\{p_i\},1)$ random variable while and add a $0$. If the outcome is a $1$ then draw from a $\mathcal{H}(I\setminus\{p_i\},0)$ random variable and add a $1$.
\end{lemma}

Again, in case all parameters $p_i$ are equal to $p$, the previous lemma can be formally expressed as

\[ P(n,p) \sim  B(1,\alpha(\{p\}_{n-1})) + B(n-1,p,1- B(1,\alpha(\{p\}_{n-1}))) . \]

Lemma \ref{evendistr} and \ref{odddistr} imply that the distributions of
$\mathcal{E}(I,\pi), \mathcal{O}(I,\pi)$ can be analyzed via the distributions
$\mathcal{H}(I\setminus \{p_i\},0)$ and $\mathcal{H}(I\setminus \{p_i\},1)$.
The next result can be used in case one is  interested in adding independent copies of $\mathcal{E}(\cdot,\cdot)$ and $\mathcal{O}(\cdot,\cdot)$.

\begin{lemma}
\label{evenpartition}
Let $I=\{p_1,\ldots,p_n\}$ and consider a partition of $I$ into disjoint, non-empty sets $I_1,I_2$.
Then the distribution of $\mathcal{H}(I,0)$ is a mixture of the independent sums $\mathcal{H}(I_1,0) + \mathcal{H}(I_2,0)$ and
$\mathcal{H}(I_1,1) + \mathcal{H}(I_2,1)$. More precisely, for an even $k$,
we have
\begin{eqnarray*}\mathbb{P}[\mathcal{H}(I,0)=k] &=& \frac{\alpha(I_1)\cdot \alpha(I_2)}{\alpha(I)} \mathbb{P}[\mathcal{H}(I_1,0) + \mathcal{H}(I_2,0)=k] \\
&+& \frac{\beta(I_1)\cdot\beta(I_2)}{\alpha(I)} \mathbb{P}[\mathcal{H}(I_1,1) + \mathcal{H}(I_2,1)=k] .
\end{eqnarray*}
\end{lemma}
\begin{proof} Write $\mathbb{P}[\mathcal{H}(I,0)=k]= \frac{\mathbb{P}[\mathcal{H}(I)=k]}{\alpha(I)}$
and note that if we regard $\mathcal{H}(I)$ as an independent sum of $\mathcal{H}(I_1)$
and $\mathcal{H}(I_2)$, then $\mathbb{P}[\mathcal{H}(I)=k]$ equals
\[ \sum_{i: i\;\text{even}} \mathbb{P}[\mathcal{H}(I_1)=i] \cdot \mathbb{P}[\mathcal{H}(I_2) =k-i] +  \sum_{i: i\;\text{odd}} \mathbb{P}[\mathcal{H}(I_1)=i] \cdot \mathbb{P}[\mathcal{H}(I_2) =k-i] . \]
Multiply and divide the  sum that runs over even indices by $\alpha(I_1)\cdot \alpha(I_2)$ and the sum that
runs over odd indices by $\beta(I_1)\cdot\beta(I_2)$ to get the result.
\end{proof}

Similarly, one can prove the following.

\begin{lemma}
\label{oddpartition}
Let $I=\{p_1,\ldots,p_n\}$ and consider a partition of $I$ into disjoint, non-empty sets $I_1,I_2$.
Then the distribution of $\mathcal{H}(I,1)$ is a mixture of the independent sums $\mathcal{H}(I_1,1) + \mathcal{H}(I_2,0)$ and
$\mathcal{H}(I_1,0) + \mathcal{H}(I_2,1)$. More precisely, for an odd $k$,
we have
\begin{eqnarray*} \mathbb{P}[\mathcal{H}(I,1)=k]&=&\frac{\alpha(I_1)\cdot \beta(I_2)}{\beta(I)} \mathbb{P}[\mathcal{H}(I_1,0) + \mathcal{H}(I_2,1)=k] \\
&+& \frac{\beta(I_1)\cdot\alpha(I_2)}{\beta(I)} \mathbb{P}[\mathcal{H}(I_1,1) + \mathcal{H}(I_2,0)=k] .
\end{eqnarray*}
\end{lemma}

The last two lemmata can be iterated. By doing so one gets that every $\mathcal{H}(I,0)$ or
$\mathcal{H}(I,1)$ random variable is a mixture
of independent sums consisting only of summands of the form $\mathcal{H}(\{a,b\},0)$, $\mathcal{H}(\{c,d\},1)$, $\mathcal{H}(\{e,f,g\},0)$ and $\mathcal{H}(\{k,l,m\},1)$, where $a,b,c,d,e,f,g,k,l,m \in (0,1)$. That is, one may apply the last two lemmata by partitioning $I$ into $I_1\cup D_1$, where $D_1$ is a doubleton. Then apply the lemma again by partitioning $I_1$ into $I_2\cup D_2$, for some doubleton $D_2$ and so on.  \\
The reason to partition $I$ this way is the next result that says that all terms of the previous mixture are
rescaled biased coins. Its proof is immediate.

\begin{lemma}
\label{rescaledcoins}
Let $I=\{p_1,p_2\}$ and $J=\{q_1,q_2,q_3\}$. Then
$\mathcal{H}(I,0)\sim 2\cdot B(1, \frac{p_1\cdot p_2}{\alpha(I)})$,
$\mathcal{H}(I,1)\sim  B(1,1)$,
$\mathcal{H}(J,0)\sim 2\cdot B(1,1- \frac{(1-q_1)\cdot(1- q_2)\cdot(1-q_3)}{\alpha(J)})$ and $\mathcal{H}(J,1)\sim 1 + 2\cdot B(1,\frac{q_1\cdot q_2\cdot q_3}{\beta(J)})$
\end{lemma}

The next result is an inequality on conditional binomial random variables.
Set $\alpha_n = \mathbb{P}[B(n,p) \; \text{even}]$ and $\beta_n=1-\alpha_n$.

\begin{lemma}
\label{conditionalbinomials} Fix a positive integer $n$ and a real number $p\in (0,1)$. Then
\[ \mathbb{P}[B(n,p,1)\geq 2k-1] \geq \mathbb{P}[B(n,p,0)\geq 2k] \]
and
\[ \mathbb{P}[B(n,p,0)\geq 2k] \geq \mathbb{P}[B(n,p,1)\geq 2k+1] . \]
\end{lemma}
\begin{proof} We induct on $n$. For $n=2$
it is easy to check that both inequalities hold true, so suppose that  both inequalities
hold true for all positive integers that are $\leq n-1$.
Let $q=1-p$. The fact that $1-2q=-1+2p$ and the symmetry of the binomial distribution imply that
it is enough to check
the inequalities for $p \in(0,1/2]$.
In order to simplify notation, set $X_{n}=B(n,p,0)$ and $Y_n=B(n,p,1)$.
From Lemma \ref{evenpartition} and Lemma \ref{oddpartition} we know that
\[ \mathbb{P}[Y_n \geq 2i-1] = \frac{p\alpha_{n-1}}{\beta_{n}}\mathbb{P}[X_{n-1}\geq 2i-2]+ \frac{(1-p)\beta_{n-1}}{\beta_n}
\mathbb{P}[Y_{n-1}\geq 2i-1] , \]
\[  \mathbb{P}[X_n\geq 2i] = \frac{(1-p)\alpha_{n-1}}{\alpha_{n}}\mathbb{P}[X_{n-1}\geq 2i]+ \frac{p\beta_{n-1}}{\alpha_n}
\mathbb{P}[Y_{n-1}\geq 2i-1] . \]
and that
\[ \mathbb{P}[Y_n\geq 2i+1] = \frac{p\alpha_{n-1}}{\beta_{n}}\mathbb{P}[X_{n-1}\geq 2i]+ \frac{(1-p)\beta_{n-1}}{\beta_n}
\mathbb{P}[Y_{n-1}\geq 2i+1] . \]
Since $p\leq 1/2$ it is easy to check that
\[ \frac{p}{\beta_n}\leq \frac{1-p}{\alpha_n} \quad \text{and} \quad \frac{p}{\alpha_n}\leq \frac{1-p}{\beta_n} . \]
Hence
\[ \mathbb{P}[Y_n\geq 2i-1]\geq \mathbb{P}[X_n\geq 2i] \]
if and only if
\begin{eqnarray*} \mathbb{P}[Y_{n-1}\geq 2i-1] \cdot\beta_{n-1}\cdot \left(\frac{1-p}{\beta_n} -\frac{p}{\alpha_n}\right) \geq\\
\mathbb{P}[X_{n-1}\geq 2i]\cdot \alpha_{n-1}\cdot
\left(\frac{1-p}{\alpha_n}- \frac{p}{\beta_n}\right)
- \frac{p \alpha_{n-1}}{\beta_n}\cdot\mathbb{P}X_{n-1}=2i-2]
\end{eqnarray*}
As $0\leq p\leq 1/2$, elementary calculations and the fact that $\alpha_n=p+(1-2p)\alpha_{n-1}$ imply
\[ \beta_{n-1}\cdot \left(\frac{1-p}{\beta_n} -\frac{p}{\alpha_n}\right) = \alpha_{n-1}\cdot
\left(\frac{1-p}{\alpha_n}- \frac{p}{\beta_n}\right) \]
and the result follows from the inductional hypothesis.
Similarly,
\[ \mathbb{P}[X_n\geq 2i] \geq \mathbb{P}[Y_n\geq 2i+1] \]
if and only if
\begin{eqnarray*} \mathbb{P}[X_{n-1}\geq 2i]\cdot \alpha_{n-1} \cdot\left(\frac{1-p}{\alpha_n}- \frac{p}{\beta_n}\right) \geq \\
\mathbb{P}[Y_{n-1}\geq 2i+1]\cdot \beta_{n-1}\cdot
\left(\frac{1-p}{\beta_n}- \frac{p}{\alpha_n}\right)
- \frac{p \beta_{n-1}}{\alpha_n}\cdot\mathbb{P}[Y_{n-1}=2i-1] .
\end{eqnarray*}
As $0\leq p\leq 1/2$, elementary calculations and the fact that $\alpha_n=p+(1-2p)\alpha_{n-1}$ imply
\[ \alpha_{n-1}\cdot \left(\frac{1-p}{\alpha_n} -\frac{p}{\beta_n}\right) = \beta_{n-1}\cdot
\left(\frac{1-p}{\beta_n}- \frac{p}{\alpha_n}\right) \]
and, once again, the inductional hypothesis finishes the proof.
\end{proof}

As a corollary we obtain the following result that will be used in our
analysis of colored coin tosses.
Recall (see \cite{SS}) that a random variable $X$ is said to be
stochastically larger than another random variable $Y$, denoted by $X\geq_{st} Y$,
if $\mathbb{P}[X \geq t] \geq \mathbb{P}[Y \geq t]$, for all $t$.

\begin{cor}\label{evenoddineqbin} Let $p_1\geq p_2\geq p$ be three real number from $(0,1)$
and fix a positive integer $n$. Then
\[ B(1,p_1)+ B(n,p, B(1,p_1))\geq_{st} B(1,p_2)+ B(n,p, B(1,p_2)) \]
and
\[ B(1,p_1)+ B(n,p, 1-B(1,p_1))\geq_{st} B(1,p_2)+ B(n,p,1- B(1,p_2)). \]
\end{cor}
\begin{proof} We only prove the first inequality, the other can be proved similarly.
Set $X_1= B(1,p_1)+ B(m,p, B(1,p_1))$ and $X_2=B(1,p_2)+ B(m,p, B(1,p_2))$. We want to prove that,
for every even integer, say $2k$, in $\{0,1,\ldots,n\}$, we have $\mathbb{P}[X_1\geq 2k]\geq \mathbb{P}[X_2\geq 2k]$.
This inequality is equivalent to
\begin{eqnarray*} p_1\cdot\mathbb{P}[B(n,p,1)\geq 2k-1]+(1-p_1)\cdot\mathbb{P}[B(n,p,0)\geq 2k]\geq \\ p_2\cdot\mathbb{P}[B(n,p,1)\geq 2k-1]+(1-p_2)\cdot\mathbb{P}[B(n,p,0)\geq 2k]
\end{eqnarray*}
and the later holds true if and only if
\[ \mathbb{P}[B(n,p,1)\geq 2k-1] \geq \mathbb{P}[B(n,p,0)\geq 2k] . \]
Lemma \ref{conditionalbinomials} finishes the proof.
\end{proof}

The following result gives a lower on a median of the random variables $A(n,p)$ and $P(n,p)$.
Recall that a \emph{median} of a random variable, $Y$,
is any number $\mu$ satisfying $\mathbb{P}[Y\geq \mu]\geq 1/2$ and
$\mathbb{P}[Y\leq \mu]\geq 1/2$. Notice that this $\mu$ might not be unique.
By abuse of notation, we will denote any median of $Y$
by $\text{Med}(Y)$.

\begin{lemma}
\label{medianevenbin} Fix a $p\in (0,1)$ and a positive integer $n$. Then a median
of a $A(n,p)$ random variable is $\geq (n-1)p-1$. Similarly, a median of a $P(n,p)$ random variable is $\geq (n-1)p-1$.
\end{lemma}
\begin{proof} We prove the result for $A(n,p)$. A similar argument works for $P(n,p)$.
For any even $k$, we have
\[ \mathbb{P}[A(n,p)\geq k] = \mathbb{P}[B(n-1,p)\geq k-1] . \]
Now it is well known (see \cite{Kaas}) that
a median of a $B(n-1,p)$ random variable is $\geq \lfloor(n-1)p\rfloor$.
If $\lfloor(n-1)p\rfloor$ is odd, then a median of $A(n,p)$ is $ \geq \lfloor(n-1)p\rfloor+1 \geq (n-1)p$.
If $\lfloor(n-1)p\rfloor$ is even, then $a:=\lfloor(n-1)p\rfloor -1$ is odd and is
such that $\mathbb{P}[B(n-1,p)\geq a]\geq 1/2$. Thus a median of $A(n,p)$ is $\geq \lfloor(n-1)p\rfloor \geq (n-1)p-1$.
\end{proof}

We will also need the following result on Bernoulli trials.

\begin{lemma}
\label{bernoullilarger} Let $p\in (0,1)$ and suppose that $X_i,i=1,\ldots,s$ are $\{0,1\}$-valued random variables
such that $\mathbb{P}[X_1=1]\geq p$ and
\begin{equation}\label{star} \mathbb{P}[X_i=1|X_1,\ldots,X_{i-1}] \geq p, \; \text{for all}\; i=2,\ldots,s .
\end{equation}
Then $\Sigma_s:=X_1+\cdots +X_s$ is stochastically larger than a $B(s,p)$ random variable. Furthermore,
it is possible to define random vectors $\mathbf{U}=(U_1,\ldots,U_s)$ and $\mathbf{V}=(V_1,\ldots,V_s)$
on a common probability space
so that the law of $(U_1,\ldots,U_s)$ is the same as the law of $(X_1,\ldots,X_s)$, each coordinate of $\mathbf{V}$ is
an independent $\text{Ber}(p)$ random variable and
\[V_i\leq U_i,\; \text{for all}\; i=1,\ldots,s,\; \text{with probability} \;  1 .\]
\end{lemma}
\begin{proof} We want to prove that
\[ \mathbb{P}[\Sigma_s \geq t] \geq \mathbb{P}[B(s,p)\geq t], \; \text{for all}\; t\in \{0,1,\ldots,s\} . \]
Note that every outcome of the random variables $X_i,i=1,\ldots,s$ is an $s$-tuple $(x_1,\ldots,x_s)\in \{0,1\}^s$.
We associate a binary vector $\mathbf{b}=(b_1,\ldots,b_s)$ to every outcome of $X_i,i=1,\ldots,s$
in such a way that the number of $1$'s in $\mathbf{b}$ has the same distribution as a $B(s,p)$ random variable. \\
To do so,  begin by drawing from $X_1$. Let $q_1=\mathbb{P}[X_1=1]$.
If $X_1=0$, then set $b_1=0$. If $X_1=1$, then let $b_1$ be the outcome of a $0/1$ coin that
shows up $1$ with probability $\frac{p}{q_1}$.
Note that $b_1=1$ with probability $p$. Now, for $i=2,\ldots,s$ do the following:
Suppose that we have sampled from $X_1,...,X_{i-1}$ and thus have formed an $(i-1)$-tuple $(x_1,...,x_{i-1})$.
Let $q_i= \mathbb{P}[X_i=1 | X_1=x_1,...,X_{i-1}=x_{i-1}]\geq p$ and now sample from $X_i$.
If $X_i=0$, then set $b_i = 0$. If $X_i=1$, then let $b_i$ be the outcome of a $0/1$ coin that
 shows up $1$ with probability $\frac{p}{q_i}$.
Notice again that $b_i =1$ with probability $p$ and this does not depend on the previous values $b_1,...,b_{i-1}$, by (\ref{star}). Thus the number of $1$'s in the vector $b=(b_1,\ldots,b_s)$ is binomially distributed.
If the vector $b$ has more than $t$ $1$'s, then also the vector $(X_1,...,X_n)$ has more than $t$ $1$'s and
first statement of the lemma follows. As $x_i\geq b_i$, for all $i=1,\ldots,s$, the second statement is immediate.
\end{proof}

The next result can be proved in a similar way.

\begin{lemma}
\label{bernoullilargertwo} Let $p\in (0,1)$ and suppose that $X_i,i=1,\ldots,s$ are $\{0,1\}$-valued random variables
such that $\mathbb{P}[X_1=1]\leq p$ and
\begin{equation}\label{star} \mathbb{P}[X_i=1|X_1,\ldots,X_{i-1}] \leq p, \; \text{for all}\; i=2,\ldots,s .
\end{equation}
Then $\Sigma_s:=X_1+\cdots +X_s$ is stochastically smaller than a $B(s,p)$ random variable. Furthermore,
it is possible to define random vectors $\mathbf{U}=(U_1,\ldots,U_s)$ and $\mathbf{V}=(V_1,\ldots,V_s)$
on a common probability space
so that the law of $(U_1,\ldots,U_s)$ is the same as the law of $(X_1,\ldots,X_s)$, each coordinate of $\mathbf{V}$ is
an independent $\text{Ber}(p)$ random variable and
\[V_i\geq U_i,\; \text{for all}\; i=1,\ldots,s,\; \text{with probability} \;  1 .\]
\end{lemma}

We end with an important result, obtained by Hoeffding (see \cite{Hoeffding}),
that will be used in the next section.

\begin{thm}\label{hoeffding} If $I=\{p_1,\ldots,p_n\}$ is a set of parameters in $(0,1)$, then
\[ \mathbb{P}[b\leq \mathcal{H}(I) \leq c] \geq \mathbb{P}[b\leq B(n,\bar{p})\leq c] , \; \text{when}\; 0\leq b \leq n\bar{p} \leq c \leq n , \]
where $\bar{p}=\frac{1}{n}\sum_{i=1}^{n}p_i$.
\end{thm}

\section{Randomly oriented graphs}\label{randorient}

Suppose that you have $n$ colors and $n$ \emph{biased} coins, all coins having the same bias.
Suppose that you color the coins in such a way that no coin has the same color on both sides.
In this section we present a method to obtain upper bounds on the median
of the number of different colors after a toss.
Note that for every such coloring of the coins
one can associate a \emph{graph}  whose
vertices correspond to the colors and whose
edges correspond to the coins. More explicitly,
for each color put a vertex in the graph and join
two vertices if and only if they are
sides of the same coin. Note  that the graph is loop-less and
that it might have parallel edges,
because the same colored coin may occur more than one time.
In addition, note that the graph may \emph{not} be connected and
that there is a
one-to-one correspondence between array of coins and
graphs and so one can choose not to distinguish
between vertices and colors as well as
between coins and edges. We call this graph the \emph{dependency} graph of the set of coins.
Fix $n$ \emph{biased} coins that are colored with $n$ colors. Let $p\in (0,1)$ be the bias of the coins and
let $G=(V,E)$ be the dependency graph of the colored coins.  Without loss of
generality we may assume that $0<p\leq \frac{1}{2}$.
Note that $|V|=|E|=n$.
Every toss of the coins gives rise to an orientation on the edges of $G$.
As a consequence, if $X_G$ is the number of different colors after the toss, then
$X_{G} = j$ corresponds to the fact that $j$ vertices in $G$ have positive in-degree, which in turn means that
$n-j$ vertices
must have in-degree $0$.  Note that none of the vertices of zero in-degree can be adjacent.
Hence if $Z_G$ is the number of vertices of zero in-degree after a toss of the coins then $X_G = n - Z_G$.
In this section we present a method to obtain an upper bound on a median of $X_G$.

In order to make an educated guess on a bound of $\text{Med}(X_G)$, one might first try to maximize $\mathbb{E}[X_G]$.
To do so, we need some extra notation.
For every vertex $v$ from $G$,  let $P_v$ be the set of edges incident to
$v$ that are oriented towards $v$ with probability $p$.
Denote  by $Q_v$ the set of edges incident to $v$ that are oriented towards
$v$ with probability $q:=1-p$. Set  $x_v=|P_v|$ and $y_v=|Q_v|$ so that $x_v+y_v= \text{deg}(v)$.

\begin{lemma} The maximum value of $\mathbb{E}[X_G]$ is  $n(1-p+p^2)$. This value is achieved by a
set of coins that uses every color twice and every color in this set appear exactly once in a $p$-side of a coin and exactly once in
a $q$-side of some other coin.
\end{lemma}
\begin{proof} Fix a graph $G$ and for every $v\in G$
denote by $C_v$ the event that vertex $v$ gets positive in-degree after a toss.
Then
\[ \mathbb{E}[X_G]= \sum_{v\in G}\mathbb{P}[C_v]=\sum_{v\in G} (1-(1-p)^{x_v}p^{y_v}) . \]
The arithmetic-geometric mean inequality implies that
\[ \sum_{v\in G} (1-p)^{x_v}p^{y_v} \geq n\cdot  (\prod_{v\in G} (1-p)^{x_v}p^{y_v})^{1/n} =  np(1-p) ,\]
since $\sum_v x_v = \sum_v y_v =n$.
We conclude that $\mathbb{E}[X_G]\leq n-np(1-p) = n(1-p+p^2)$. The second statement is immediate.
\end{proof}

Notice that the graph $G$ for which the mean of $X_G$ is maximum is a union of cycles.
Note also that the function $f(p)=1-p+p^2, p\in (0,1)$
is convex and attains its minimum at $p=\frac{1}{2}$.
This means that the maximum mean is minimized when $p=\frac{1}{2}$.

Finding an upper bound on a median of $X_G$ turns out to be more involved.
Our main result on the median of $X_G$ is the following.

\begin{thm}
\label{medianbound}
For any loop-less multi-graph $G$ on $n$ vertices and $n$ edges,
a median of $X_G$ is $\leq n- \frac{p^2}{1+(1-2p)^2}n + \frac{3}{4}$.
\end{thm}

The rest of the section is devoted to the proof of this theorem. We will analyze the distribution of $X_G$
via the distribution of $E_G$, the number of vertices with even in-degree after a toss of the coins.
The reason to do so is contained in the following result.

\begin{lemma}
\label{evenindegreetwo} Fix a (possibly disconnected) graph $G$, on $n$ vertices and $n$ edges
as well as an  orientation on the edges of $G$. Let $Z_G$ be the number of vertices of zero
in-degree and $E_G$ the number of vertices of even in-degree in $G$. Then
\[ Z_G \geq \frac{1}{2} E_G . \]
A lower bound on $\text{Med}(E_G)$ gives an upper bound on $\text{Med}(X_G)$. More precisely,
\[ \text{Med}(X_G) \leq n - \frac{1}{2} \text{Med}(E_G) . \]
\end{lemma}
\begin{proof} Let $Y_G=E_G-Z_G$. For $i=1,2,\ldots,$ set $I_i:=\{v\in G: \text{deg}^{-}(v)=i\}$.
From the in-degree sum formula we have that
\[ n =    \sum_{v\in G}    \text{deg}^{-}(v) = \sum_{i\geq 1} i |I_i| .     \]
In addition, $n=Z_G+ \sum_{i\geq 1}|I_i|$. Hence
\begin{eqnarray*} n-n &=&  \sum_{i\geq 1} i|I_i| - \sum_{i\geq 1}|I_i| - Z_G \\
                      &=& \sum_{i\geq 1}(i-1)|I_i|  - Z_G \\
                      &\geq&  Y_G - Z_G \\
                      &=& E_G-2Z_G ,
\end{eqnarray*}
which implies that $2Z_G \geq E_G$, thus proving the first statement. From this we can conclude that
\[ X_G = n-Z_G \leq n-\frac{1}{2}E_G, \]
and so $\text{Med}(X_G) \leq n - \frac{1}{2} \text{Med}(E_G)$, as required.
\end{proof}

The idea behind looking at the number of vertices of even in-degree is the following.
Recall that we are interested in obtaining an upper bound on a median of $X_G$.
Since $X_G=n-Z_G$, the problem is equivalent to
obtaining a lower bound on a median of $Z_G$.
From the previous lemma we know that $Z_G \geq \frac{1}{2}E_G$, for all oriented graphs $G$.
This means that if we can determine a lower bound on a median of $E_G$ then we will also
have obtained an upper bound on a median of $X_G$, by Lemma \ref{evenindegreetwo}.
Furthermore, in case $G$ is connected, one can "estimate" the distribution of $E_G$ from below.
More precisely, let  $E_G$ be the number of vertices with even in-degree after a random orientation on the edges of $G$.
Recall (see \cite{SS}) that if $X$ and $Y$ are random variables, then we say that
$X$ is \emph{stochastically larger} than $Y$, denoted by $X\geq_{st} Y$,  if
\[\mathbb{P}[X\geq t]\geq \mathbb{P}[Y\geq t], \; \text{for all}\; t . \]
In case $\mathbb{P}[X\geq t]= \mathbb{P}[Y\geq t]$, for all $t$ we will write $X=_{st} Y$.
Our main result on the distribution of $E_G$ is the following.

\begin{thm}
\label{distrevenindegreenew} Suppose that $G=(V,E)$ is a connected multi-graph on $n$ vertices
and $m\geq n-1$ edges. Let $d_v$ be the degree of vertex $v$, set
$\pi_v:=\frac{d_v}{2m}$ and let $\pi$ be the probability vector with coordinates $\pi_v, v\in V$.
Let $E_{G}$ be the number of
even in-degree vertices after orienting each edge towards its tail with probability $p$ and towards
its head with probability $1-p$. Assume $p<1-p$ and let $\{p\}_n$ be the set consisting of $n$ copies of $p$.
Then, if $m-n$ is even, $E_{G}$ is stochastically larger than
a $\mathcal{E}(\{p\}_n,\pi)$ random
variable. If $m-n$ is odd, then $E_{G}$ is stochastically larger
than a $\mathcal{O}(\{p\}_n,\pi)$ random
variable.
\end{thm}

Note that by the remarks following the definition of even-sum (resp. odd-sum) toss of $n$ coins, we know that
\[ \mathcal{E}(\{p\}_n,\pi)\sim A(n,p) \quad \text{and} \quad \mathcal{O}(\{p\}_n,\pi)\sim P(n,p) . \]

We prove this Theorem in a series of lemmata. We begin with
a result that imposes parity restrictions on $E_G$.
Denote by $\text{deg}^{-}(v)$ the in-degree of vertex $v$.

\begin{lemma}
\label{degreesum}
Suppose that $G$ is a (possibly disconnected) graph on $n$ vertices
and $m$ edges. Fix some orientation on the edges and
let $O_{G}, E_{G}$ be the number of odd and even in-degree vertices respectively.
Then the parity of $E_{G}$ equals the parity of $m-n$.
\end{lemma}
\begin{proof}
The in-degree sum formula states that
\[ \sum_{v\in G} \text{deg}^{-}(v) = m . \]
From this we get that the parity of $O_{G}$ equals the parity of $m$.
As $n-E_{G}= O_{G}$, it follows that the parity of $m$ equals the parity of $n-E_{G}$,
as required.
\end{proof}

The following labeling on the vertices
and edges of a tree will also be of use.
Recall that a \emph{leaf} in a tree is a vertex of degree $1$.

\begin{lemma}
\label{treeenumerationtwo} Let $T$ be a tree on $n$ vertices and fix any edge $f\in T$.
Then there exists a labeling, $v_1,\ldots ,v_n$, of the vertices
and a labeling, $e_1,\ldots, e_{n-1}$, of the edges of $T$
such that\\
(i) edge $f$ has label $e_{n-1}$; \\
(ii) the only edge incident to
vertex $v_i, i=1,\ldots , n-1$, among the edges with labels $\{e_i,e_{i+1}, \ldots , e_{n-1}\}$ is
the edge with label $e_i$.
\end{lemma}
\begin{proof} The statement is clearly true if $n=2$, so suppose that $n>2$.
Fix a tree, $T$, on $n>2$ vertices and choose any of its edges. Label this edge $e_{n-1}$
and label its endpoints  $v_n$ and $v_{n-1}$ arbitrarily.
Notice that not both $v_n$ and $v_{n-1}$ can be leaves.
If $v_n$ or $v_{n-1}$ is a leaf, say $v_{n}$, then
consider the vertex set $L$ of
leaves in $T$ except $v_n$ and label them $v_1, v_2,\ldots ,v_{\ell}$.
If $v_n$ is not a leaf, then consider all leaves of $T$ and label them in the same manner.
Note that $L$ is not empty even if $v_n$ is a leaf since any tree with at least two vertices has at least
two leaves.
Now label each edge incident to $v_j$ with $e_j$, for $j=1,2,\ldots ,\ell$.
Now consider the tree $T':=T\setminus \{v_1, v_2,\ldots ,v_{\ell}\}$ and
repeat this process on the leaves of $T'$ again sparing $v_n$ or $v_{n-1}$ if it is a leaf of $T'$.
We keep on labeling the leaves and edges of
the subtrees until we end up with the graph consisting of the edge $e_{n-1}$ only.
It is evident that the labeling satisfies the required condition.
\end{proof}

Note that we can label any edge of $T$ with $e_{n-1}$
and any endpoint of $e_{n-1}$ with $v_n$.
We will call a labeling on the vertices and edges of a tree, a \emph{good labeling}
if it satisfies the conditions of Lemma \ref{treeenumerationtwo}. Notice also that
if we are given a good labeling of a tree and we interchange the labels
$v_n$ and $v_{n-1}$ then we get another good labeling of the same tree.
We collect this observation in the following.

\begin{lemma}
\label{goodenumer} Let $T$ be a tree on $n$ vertices and fix two
adjacent vertices $u_1,u_2$ of $T$. Suppose that
$T$ has a good labeling such that $u_1$ has label $v_{n-1}$ and $u_2$ has label $v_n$.
Then the labeling that interchanges the labels of $u_1$ and $u_2$ and keep all other labels the same is also
a good labeling.
\end{lemma}

Note that the previous lemma says that for any edge $f=(u,w)$ of $T$ there is
a one-to-one correspondence between good labelings for which $u$ gets the label $v_n$ and $w$ gets label $v_{n-1}$
and good labelings for which $u$ gets the label $v_{n-1}$ and $w$ gets label $v_n$.
We will also need the following observation
on the spanning trees of connected graphs.

\begin{lemma}
\label{edgetree} Suppose that $G=(V,E)$ is a connected graph and fix any edge $e\in E$.
Then there exists a spanning tree, $T$, of $G$  such that $e$ is an edge of $T$, i.e. $e\in T$.
\end{lemma}
\begin{proof} Let $T=(V,E')$ be a spanning tree of $G$. If $e\in E'$ then we are done, so
suppose that $e\notin E'$. This means that if we add $e$ to $E'$ then we create a cycle.
Now note that if we delete any edge, $e'\neq e$, from this cycle we get a spanning
tree $T'$ of $G$ for which $e$ belongs to $T'$.
\end{proof}

After each assignment of orientation to the edges, let $x^{-}_{v}$ be  the number of edges in $P_v$ that are
oriented towards $v$, and  $y_{v}^{-}$ be the number of edges in
$Q_v$ that are oriented towards $v$.
In the following result we compute the probability that a certain vertex has even in-degree.

\begin{lemma}
\label{probevenindeg} If $v\in V$ is such that $y_v$ is even, then
\[ \mathbb{P}[\text{deg}^{-}(v) \; \text{even}] =  \mathbb{P}[B(\text{deg}(v),p) \; \text{even}]. \]
If $v\in V$ is such that $y_v$ is odd, then
\[  \mathbb{P}[\text{deg}^{-}(v) \; \text{even}] =  \mathbb{P}[B(\text{deg}(v),p) \; \text{odd}]
\]
\end{lemma}
\begin{proof} We only prove the first equality. The second can be proved similarly.
Note that $\text{deg}^{-}(v)$ is even if and only if either both $x^{-}_{v}$ and $y_{v}^{-}$ are even, or
both are odd. Thus
\begin{eqnarray*}
\mathbb{P}[\text{deg}^{-}(v) \; \text{even}] &=&  \mathbb{P}[x^{-}_{v} \; \text{even}] \cdot \mathbb{P}[y^{-}_{v} \; \text{even}]\\
&+& \mathbb{P}[x^{-}_{v} \; \text{odd}] \cdot \mathbb{P}[y^{-}_{v} \; \text{odd}] ,
\end{eqnarray*}
which can be rewritten as
\begin{eqnarray*}
\mathbb{P}[\text{deg}^{-}(v) \; \text{even}]&=&   \mathbb{P}[B(x_v,p)\; \text{even}] \cdot \mathbb{P}[B(y_v,1-p)\; \text{even}] \\
&+& \mathbb{P}[B(x_v,p)\; \text{odd}]\cdot \mathbb{P}[B(y_v,1-p)\; \text{odd}]
\end{eqnarray*}
and so $\mathbb{P}[\text{deg}^{-}(v) \; \text{even}]$ equals
\[
\frac{1}{2}(1+(1-2p)^{x_v}) \cdot \frac{1}{2}(1+(1-2q)^{y_v}) + \frac{1}{2}(1-(1-2p)^{x_v}) \cdot \frac{1}{2}(1-(1-2q)^{y_v}) . \]
Now from the fact that $1-2q=-1+2p$ and $y_v$ is even, we can conclude that the last expression is the same as
\[ \frac{1}{2}(1+(1-2p)^{x_v})\cdot \frac{1}{2}(1+(1-2p)^{y_v})+ \frac{1}{2}(1-(1-2p)^{x_v}) \cdot \frac{1}{2}(1-(1-2p)^{y_v}) \]
which in turn is equal to
\[ \frac{1}{2} + \frac{1}{2} (1-2p)^{\text{deg}(v)}  \]
and proves the lemma.
\end{proof}

The next result is crucial since it will reduce the problem of obtaining an upper bound on a median of $X_G$
to the one of obtaining a lower bound on a median of a conditional binomial distribution.
Recall that we assume $p\leq 1/2$.

\begin{lemma}
\label{independentparitynew}
Fix some vertex $v$ of the graph, fix an edge, $e$, that is incident to $v$ and
let $C$ be the set consisting of all edges edges incident to $v$ except $e$.
Let $C^{-}$ denote the number of edges from $C$ that are oriented towards $v$ after a toss.
Then
\[ \mathbb{P}[\text{deg}^{-}(v)\;\text{even}|C^{-}] \geq p . \]
\end{lemma}
\begin{proof}
Suppose the coins corresponding to $C$ have
been flipped.
Let $C^{-}$ be the number of edges in $C$
which are oriented towards $v$ after the toss. Suppose that the edge $e$ corresponds to a coin that
is oriented towards $v$ with probability $p$. The other case is similar. Then
\begin{eqnarray*} \mathbb{P}[\text{deg}^{-}(v) \; \text{even}| C^{-}] &=& (1-p)\cdot \mathbf{1}_{\{C^{-}\; \text{even}\}} + p\cdot \mathbf{1}_{\{C^{-}\; \text{odd}\}}  \\
&=& p + (1-2p)\cdot\mathbf{1}_{\{C^{-}\; \text{even}\}} \\
&\geq& p,
\end{eqnarray*}
where $\mathbf{1}_{\{\cdot\}}$ denotes indicator.
Note that in case $p=\frac{1}{2}$ the last inequality is in fact equality and that the same
computation shows that $\mathbb{P}[\text{deg}^{-}(v)\;\text{even}|C^{-}] \leq 1-p$.
\end{proof}

For every vertex $v\in G$, denote by $\theta_v$ the probability that the in-degree
of $v$ is even. Note that, by  Lemma \ref{probevenindeg}, $\theta_v$ is either equal to
$\mathbb{P}[\text{Bin}(\text{deg}(v),p)\; \text{even}]$ or to $\mathbb{P}[\text{Bin}(\text{deg}(v),p)\; \text{odd}]$.
Thus $p\leq \theta_v\leq 1-p$, for all $v\in V$. \\

We now have all the necessary tools to prove our main result on $E_G$.

\begin{proof}[Proof of Theorem \ref{distrevenindegreenew}]
Recall that for every edge we toss a coin to decide on its orientation.
All these $m$ coins, $c_i,i=1\ldots,m$, are independent. Since the order with which we toss the coins doesn't matter
we may, equivalently,
suppose that we toss the coins in the following way:
we choose a coin, say coin $c_i$,
with probability $\frac{1}{m}$,
flip the remaining $m-1$ coins in any way we want and then toss the coin $c_i$.
Tossing this way does not affect the distribution of $E_G$ but allows us to use
Lemma \ref{treeenumerationtwo}. More precisely, we may suppose that once the coin $c_i$
is chosen, then we toss the remaining $m-1$ coins according to a good labeling,
$v_1,\ldots ,v_n; e_1,\ldots, e_{n-1}$, of a spanning tree $T$ of $G$
that contains the edge corresponding to $c_i$, say this edge is $f_i=[u,w]$,
and with the good labeling of $T$ chosen in such a way that the edge $f_i$ gets label $e_{n-1}$;
we can use this specific good labeling of $T$ and
first toss the coins corresponding to edges that do \emph{not} belong to $T$
in any way we like and then
toss the coins that correspond to edges $e_1,\ldots, e_{n-1}$ in that specific order. This way
the coin $c_i$ is flipped last and we do not affect the distribution of $E_G$.
Note that, by Lemma \ref{edgetree},
there exists a spanning tree, $T$, of $G$
containing edge $f_i$ and we can always construct
a good labeling of $T$ for which $f_i$ gets label $e_{n-1}$, by Lemma \ref{treeenumerationtwo}.
Furthermore, the edge $f_i$ has two endpoints, $u,w$, and the probability
that vertex $u$ has label $v_n$ equals $1/2$, by Lemma \ref{goodenumer}.
Since we fix coin $c_i$ with probability $1/m$ it follows that, for every vertex $v\in V$,
the probability that we toss the coins according to a good labeling of a spanning tree $T$ of $G$
for which
vertex $v$ gets label $v_n$ equals $\frac{d_v}{2m}$.\\
So let $T$ be a spanning tree of $G$ with a good labeling and
recall that we are going to do the following: first we randomly orient
the edges that do \emph{not} belong to $T$
and then randomly orient the edges $e_1, e_2, \ldots , e_{n-1}$ in that order.
Note that the probability that the vertex with label $v_1$ has even in-degree equals $\theta_{v_1}\geq p$.
The fact that $T$ has a good labeling implies that, for $j=1,\ldots,n-1$, once the edge $e_j$ is given
an orientation, then the
parity of vertex $v_j$ is determined.
Lemma \ref{independentparitynew} gives that
once the parity of vertex $v_j$ is determined, the probability that
vertex $v_{j+1}$ has even in-degree is $\geq p$. Only
the parity of the vertex with label $v_n$ is deterministic given the parities of the previous vertices.
Let $\delta_i$  be the indicator of the event $\{\text{deg}^{-}(v_i)\;\text{is even}\}$, for
$i=1,2,\ldots, n$. Thus $E_G=\delta_1+\cdots+\delta_n$
and each $\delta_i,i=1,\ldots,n-1$ is stochastically larger than a  $B(1,p)$
random variable. From Lemma \ref{bernoullilarger} we know that there exist random binary vectors $\mathbf{U}=(U_1,\ldots,U_{n-1})$ and $\mathbf{V}=(V_1,\ldots, V_{n-1})$ defined on a common probability
space such that the law of $\mathbf{U}$ is the same as the law of $(\delta_1,\ldots,\delta_{n-1})$,
each $V_i$ is an independent Bernoulli $\text{Ber}(p)$ random variable and
\[ \sum_{i=1}^{n-1} U_i \geq \sum_{i=1}^{n-1}V_i \quad \text{with probability}\; 1 . \]
In addition we know that $\sum_{i=1}^{n-1}V_i\sim B(n-1,p)$.
To end the proof, suppose that $m-d$ is even. The other case is similar.
Thus $E_G$ is even as well and $E_G\sim U_1 + \cdots + U_{n-1} + \delta_n$,
where $\delta_n=1$ if $U_1+\cdots + U_{n-1}$ is odd and $\delta_n=0$ if $U_1+\cdots + U_{n-1}$ is even.
Now let $\gamma_n=1$ if $V_1+\cdots+V_{n-1}$ is odd and $\gamma_n=0$ if $V_1+\cdots + V_{n-1}$ is even, in order
to guarantee that $V_1+\cdots +V_{n-1} + \gamma_n$ is always even.
Since $U_1+\cdots + U_{n-1}\geq V_1+\cdots+V_{n-1}$ with probability $1$, we also have that
$U_1+\cdots +U_{n-1}+\delta_n \geq V_1+\cdots +V_{n-1} +\gamma_n$ with probability $1$
and the result follows.
\end{proof}

Note that in case $p=\frac{1}{2}$ Lemma \ref{independentparitynew} gives that
once the parity of vertex $v_j$ is determined, the probability that
vertex $v_{j+1}$ has even in-degree is equal to $\frac{1}{2}$, and so the parity of $v_{j+1}$
is \emph{independent}  of the parity of $v_1, v_2, \ldots ,v_{j-1}$. Only
the parity of $v_n$ is deterministic given the parities of the previous vertices.
This implies that the random variables $\delta_i,i=1,\ldots, n-1$ in the proof of Theorem
\ref{distrevenindegreenew} satisfy $\delta_1+\cdots +\delta_{n-1} =_{st} B(n-1,1/2)$ and the following
result (which is Theorem $4$ in \cite{PelSch}) follows.

\begin{cor}\label{fairdistr} Suppose that $p=\frac{1}{2}$. If $m-n$ is even, then $E_G$ has the same
distribution as a $A(n,1/2)$ random variable. If $m-n$ is odd, then $E_G$ has the same
distribution as a $P(n,1/2)$ random variable.
\end{cor}

Using Lemma \ref{medianevenbin} and Lemma \ref{evenindegreetwo} we have the following
result on $X_G$, in case $G$ is connected.

\begin{cor} Let $G$ be a connected loop-less multi-graph on $n$ vertices and $n$ edges.
Then a median of $X_G$ is $\leq n-\frac{1}{2}(n-1)p+\frac{1}{2}$.
\end{cor}

We now turn to the proof of Theorem \ref{medianbound}.
Recall that the dependency graph $G=(V,E)$
of the colored coins might not be connected.
Suppose it consists of $t$
connected components, $G_1,\ldots,G_t$, each having $n_i$ vertices and $m_i$ edges such that
$\sum n_i=n$ and $\sum m_i=n$.
Let also $E_{G_i}$ be the number of vertices of even in-degree
in each component, after a toss.  Hence the total number of vertices of even in-degree after a toss, $E_G$ is
equal to the independent sum
$E_{G_1}+\cdots +E_{G_t}$. As $|V|=|E|=n$, it follow from Lemma \ref{degreesum} that $E_G$ is even.
By Theorem \ref{distrevenindegreenew}, the distribution of each $E_{G_i}$ is
stochastically larger than
a $A(\cdot,p)$ or $P(\cdot,p)$ random variable. More precisely, suppose that
the first $t_1$ components of $G$ correspond to a $A(\cdot,p)$ random variable and the remaining $t_2$
components correspond to a $P(\cdot,p)$ random variable, so that $t_1+t_2=t$ and $t_2$ is even.
 Let $\{p\}_{k}$ denote the set consisting of $k$ parameters that are all equal to $p$. From Theorem \ref{distrevenindegreenew} we know that
\[ E_{G_i} \geq_{st}  A(n_i,p),\; \text{for}\; i=1,\ldots,t_1 \]
and
\[ E_{G_i} \geq_{st}  P(n_i,p),\text{for}\;  i=t_1+1,\ldots,t .\]
Hence, the total number of even in-degree vertices, $E_G$, is   stochastically larger than the independent sum
\[  \sum_{i=1}^{t_1} A(n_i,p) + \sum_{i=t_1+1}^{t} P(n_i,p) . \]

Since  $p\in(0,1/2]$ we have
$\beta(\{p\}_{n_i-1})\geq p$ and $\alpha(\{p\}_{n_i-1})\geq p$ and thus Corollary \ref{evenoddineqbin} implies that
\begin{eqnarray*} A(n_i,p)&\sim& B(1,\beta(\{p\}_{n_i-1})) + B(n_i-1,p, B(1,\beta(\{p\}_{n_i-1})))\\
                          &\geq_{st}& B(1,p) + B(n_i-1,p, B(1,p))
\end{eqnarray*}
and
\begin{eqnarray*} P(n_i,p)&\sim& B(1,\alpha(\{p\}_{n_i-1})) + B(n_i-1,p,1- B(1,\alpha(\{p\}_{n_i-1})))\\
                          &\geq_{st}& B(1,p) + B(n_i-1,p,1- B(1,p))
\end{eqnarray*}
and so $E_G$ is stochastically larger than
\[ \sum_{i=1}^{t_1} B(1,p)+ B(n_i-1,p,B(1,p)) + \sum_{i=t_1+1}^{t} B(1,p)+ B(n_i-1,p,1-B(1,p))
\]
This independent sum  takes even values (recall $t_2$ is even) and, by Lemma \ref{evendistr} and Lemma  \ref{odddistr}, is equivalently described as follows.
Toss $t$ independent $0/1$ coins, $c_i,i=1,\ldots,t$,
each having probability $p$ of landing on $1$.
Let $\Gamma=(\gamma_1,\ldots,\gamma_t)\in \{0,1\}^t$ be a particular outcome of the toss.
This is a binary
vector of length $t$. If $B_{\Gamma}$ is the
number of $1$'s in this vector, then add
$B_{\Gamma}$ to the outcome of the independent sum
\[   \mathcal{H}|\Gamma :=  \sum_{i=1}^{t_1}B(n_i-1,p,\gamma_i) + \sum_{i=t_1+1}^{t}B(n_i-1,p,1-\gamma_i) , \]
thus forming the sum $B_{\Gamma}+\mathcal{H}|\Gamma$.
Note that $B_{\Gamma}\sim B(t,p)$.
Now each particular vector $\Gamma$ can be equivalently obtained in the following way.
Fist toss a coin with probability of success $\frac{1}{2}(1+(1-2p)^t)$. If the outcome is a success,
then arrange $t$ independent $0/1$ coins (whose probability of landing on $1$ equals $p$) on a line and toss them until you see an \emph{even} number of $1$'s.
If $\Gamma_{e}$ is the resulting binary vector and $B_e$ is
the number of $1$'s in $\Gamma_{e}$, then $B_e \sim B(t,p,0)$ and $B_{\Gamma}+\mathcal{H}|\Gamma$
equals $B_{e}+\mathcal{H}|\Gamma_{e}$ with probability $\frac{1}{2}(1+(1-2p)^t)$.
If the outcome is a failure,
then toss $t$ independent $0/1$ coins until you see an \emph{odd} number of $1$'s.
If $\Gamma_{o}$ is the resulting binary vector
and $B_o$ is the number of $1$'s in $\Gamma_{o}$, then $B_o \sim B(t,p, 1)$ and $B_{\Gamma}+\mathcal{H}|\Gamma$
equals $B_{o}+\mathcal{H}|\Gamma_{o}$ with probability $\frac{1}{2}(1-(1-2p)^t)$.
Hence $B_{\Gamma}+\mathcal{H}|\Gamma$ is a mixture of the sums $B_{e}+\mathcal{H}|\Gamma_{e}$
and $B_{o}+\mathcal{H}|\Gamma_{o}$.

\begin{lemma}\label{last} A median of $B_{\Gamma}+\mathcal{H}|\Gamma$
is $\geq n\bar{p} -\frac{3}{2}$, where $\bar{p}:= \frac{2p^2}{1+(1-2p)^2}$.
\end{lemma}
\begin{proof} First toss a coin to decide whether you take a vector, $\Gamma_{e}$, with an even
number of $1$'s or a vector, $\Gamma_{o}$, with an odd number of $1$'s.
Suppose that we end up with a vector
$\Gamma_{e}$. The other case is similar.
This vector
gives rise to the  sum $ B_{e}+\mathcal{H}|\Gamma_{e}$. Then $B_e\sim B(t,p, 0)$ and
each term in $\mathcal{H}|\Gamma_e$ is of the form $B(n_i-1,p,0)$ or
$B(n_i-1,p,1)$. Apply lemmata \ref{evenpartition} and \ref{oddpartition} repeatedly to write each term
of the sum  $B_{e}+\mathcal{H}|\Gamma_{e}$ as a mixture of independent sums consisting only
of terms $\mathcal{H}(J,0)$ and $\mathcal{H}(J,1)$ for which $|J|$ equals $2$ or $3$. Thus the initial sum, $B_{e}+\mathcal{H}|\Gamma_{e}$,
is a mixture of independent sums of terms $\mathcal{H}(J,0)$ and $\mathcal{H}(J,1)$ for which $|J|$ equals $2$ or $3$. 
To end the proof, we show that a median of any independent sum in this mixture is $\geq \geq n\bar{p}-\frac{3}{2}$. 
Suppose that $\Xi$ is a particular independent sum consisting of $a$ terms of the form $B(2,p,0)$, $b$ terms
of the form  $B(2,p,1)$, $c$ terms of the form $B(3,p,0)$ and $d$ terms of the form
$B(3,p,1)$.  Thus $2a+2b+3c+3d=n$.
Lemma \ref{rescaledcoins} implies that
\[ B(2,p,0)\sim 2\cdot B\left(1,\bar{p} \right), \; B(2,p,1)\sim B(1,1), \]
where $\bar{p}= \frac{2p^2}{1+(1-2p)^2}$, and that
\[ B(3,p,0)\sim 2\cdot B\left(1,\hat{p} \right), \; B(3,p,1)\sim 1 +2 \cdot B\left(1,\tilde{p}\right), \]
where $\hat{p}=\frac{6p^2(1-p)}{1+(1-2p)^3}$ and $\tilde{p}=\frac{2p^3}{1-(1-2p)^3}$. Denote
\[ \Psi:= B(a,\bar{p}) + B(c,\hat{p}) + B(d,\tilde{p}) . \]
Then $\Xi = 2\Psi + b+d$ and so a median of $\Xi$ can be estimated via a median of $\Psi$. Hence a median
of $\Xi$ is $\geq n\bar{p}-\frac{3}{2}$ if and only if a median of $\Psi$ is $\geq \frac{n\bar{p}-b-d}{2}-\frac{3}{4}$.
Using the fact the $2a+2b+3c+3d=n$ we can write
\[  \frac{n\bar{p}-b-d}{2} =  a\bar{p}+ b\left(\bar{p}-\frac{1}{2}\right)+c \frac{3\bar{p}}{2} + d\left(\frac{3\bar{p}}{2}-\frac{1}{2}\right):=\mu_{\ast} . \]
Note that $\bar{p}-\frac{1}{2}\leq 0$.
As $0\leq p\leq 1/2$, elementary calculations show that $\hat{p}\geq \frac{3\bar{p}}{2}$ and $\tilde{p}\geq  \frac{3\bar{p}}{2}-\frac{1}{2}$.
This implies that
\[ \mathbb{E}[\Psi]= a\bar{p}+c\hat{p}+d\tilde{p}\geq  \mu_{\ast} . \]
From Hoeffding's result (Theorem \ref{hoeffding}) we know that
\[ \mathbb{P}\left[\Psi\geq \mu_{\ast}-\frac{3}{4}\right] \geq \mathbb{P}\left[B(a+c+d,p_0)\geq \mu_{\ast}-\frac{3}{4}\right] ,  \]
where $p_0=\frac{1}{a+c+d}(a\bar{p}+c\hat{p}+d\tilde{p})$ and so it is enough to show that a median of
a $B(a+c+d,p_0)$ random variable is $\geq \mu_{\ast} -\frac{3}{4}$.
Now, it is well known (see \cite{Hamza}) that the smallest uniform
(with respect to both parameters) distance between the mean and a median of a binomial distribution is $\leq \ln 2<\frac{3}{4}$.
This means that a median of $B(a+c+d,p_0)$ is $\geq a\bar{p}+c\hat{p}+d\tilde{p}-\frac{3}{4} \geq \mu_{\ast} -\frac{3}{4}$
and the lemma follows.
\end{proof}

The proof of the main result of this section is almost complete.

\begin{proof}[Proof of Theorem \ref{medianbound}]
Since $E_G$ is stochastically larger than $B_{\Gamma} + \mathcal{H}|\Gamma$ and
a median of $B_{\Gamma} + \mathcal{H}|\Gamma$ is $\geq n\bar{p} -\frac{3}{2}$,
we conclude that the median of $E_G$ is
$\geq n\bar{p} -\frac{3}{2}$.
Theorem \ref{medianbound} follows since, from Lemma \ref{evenindegreetwo}, we have
\[ \text{Med}(X_G)\leq n-\frac{1}{2}\text{Med}(E_G)\leq n-\frac{n}{2}\bar{p}+\frac{3}{4} . \]
\end{proof}

We end this section by noticing that our method works also in case one is interested in
estimating  $X_G$ from below. Since $X_G\geq n-E_G$, for all graphs $G$ it is enough to estimate
the probability distribution of $E_G$ from above, i.e., to find a random variable
that is stochastically larger than $E_G$. We know that $\theta_v\leq 1-p$, for all $v\in V$
and a modification of the proof of Theorem \ref{distrevenindegreenew} along with Lemma
\ref{bernoullilargertwo} shows that the following is true.

\begin{thm}\label{distrevenindegreenewtwo} Suppose that $G=(V,E)$ is a connected multi-graph on $n$ vertices
and $m\geq n-1$ edges. Let $d_v$ be the degree of vertex $v$, set
$\pi_v:=\frac{d_v}{2m}$ and let $\pi$ be the probability vector with coordinates $\pi_v, v\in V$.
Let $E_{G}$ be the number of
even in-degree vertices after orienting each edge towards its tail with probability $p$ and towards
its head with probability $1-p$. Assume $p<1-p$ and let $\{1-p\}_n$ be the set consisting of $n$ copies of $1-p$.
Then, if $m-n$ is even, $E_{G}$ is stochastically smaller than
a $\mathcal{E}(\{1-p\}_n,\pi)$ random
variable. If $m-n$ is odd, then $E_{G}$ is stochastically smaller
than a $\mathcal{O}(\{1-p\}_n,\pi)$ random
variable.
\end{thm}

\section{Random graphs}\label{randgraphsection}

In this section we apply our method to the distribution of the number of vertices with
odd degree in random sub-graphs of fixed graphs. More precisely, let $G$ be any \emph{connected graph} on $n$ vertices
and for each edge of $G$ toss a coin that shows up tails with probability $p$, independently
for all edges. If the result of the coin is tails, then keep the edge. If the result is heads, delete the edge.
The distribution of the vertex degree in such models has been well studied (see \cite{Bollobas} for a whole chapter on this topic).
The resulting sub-graph of $G$ that remains after the toss of the coins is random.
Let $q=1-p$ and denote by
$O_{n,p}(G)$ the number of vertices of \emph{odd degree} in the resulting graph.
The following folds true.

\begin{thm} If $0\leq p\leq \frac{1}{2}$ then the random variable $O_{n,p}(G)$ is stochastically larger than
a $A(n,p)$ random variable. If $\frac{1}{2}\leq p \leq 1$, then $O_{n,p}(G)$ is stochastically larger than
a $A(n,q)$ random variable.
\end{thm}
\begin{proof}
The proof is similar to the proof of Theorem \ref{distrevenindegreenew}, so we only sketch it.
Suppose that $0\leq p\leq \frac{1}{2}$. The other case is similar.
Let $T$ be a spanning subgraph of $G$ with a good labeling, $v_1,\ldots,v_n;e_1,\ldots,e_{n-1}$ on
its vertices and edges given by Lemma \ref{treeenumerationtwo}.
By Lemma  \ref{binparity}
we know that the probability that $d_i:=\text{deg}(v_i)$ is odd is equal to  $\frac{1}{2}(1-(1-2p)^{d_i})$, for $i=1,\ldots,n$.
Toss all coins to decide which edges are included in the sub-graph, except the coins corresponding
to the edges $e_i,i=1,\ldots,n-1$.
Now begin from vertex $v_1$ and toss a coin to decide whether edge $e_1$ is included or not.
Then proceed to vertex $v_2$ and toss a coin to decide on the edge $e_2$, and in general,
at step $j,j=1,\ldots,n-1$ move from vertex $j-1$ to vertex $j$ and toss a coin to decide if edge $e_j$ is included or not.
Let $C_j$ be the set of edges that are included in the graph and are incident to $v_j$ at step $j-1$.
As in Lemma \ref{independentparitynew},
by conditioning on whether $|C_j|$ is even or odd we conclude that
\[ \mathbb{P}[\text{deg}(v_j)\; \text{odd}|C_j]\geq p , \]
Hence the parity of each vertex $v_j,j=1,\ldots,n-1$, is stochastically larger than a $B(1,p)$ random variable.
Only the parity of vertex $v_n$ is deterministic, given the parities of the previous vertices.
The result follows from the fact that the degree-sum formula implies that $O_{n,p}$ is even.
\end{proof}

Notice that in case $p=\frac{1}{2}$ we obtain the following result.

\begin{cor}\label{randgraph} If $p=\frac{1}{2}$ then, for any
connected graph $G$, $O_{n,1/2}(G)$ has the same distribution
as a $A(n,1/2)$ random variable.
\end{cor}

\section{Some applications}\label{apll}

Let $G=(V,E)$ be a connected undirected graph and fix $T\subseteq V$. An orientation
of $G$, is an assignment of direction to each edge of $G$. An orientation of $G$
is called $T$-odd if the vertices in $T$ are the only ones having odd in-degree.
We allow $T$ to be the empty set in
which case $\emptyset$-odd orientation simply means that all vertices of $G$ have even in-degree.
The following result is obtained in \cite{Frank}, using induction.

\begin{lemma} A connected graph, $G=(V,E)$,
on $n$ vertices and $m$ edges has a $T$-odd orientation if and only if $|T|+|E|$ is even.
\end{lemma}
\begin{proof} Suppose first that $G$ has a $T$-odd orientation.  Let $E_G$ be the number
of even in-degree vertices, $O_G$ the number of odd in-degree vertices. From Lemma \ref{degreesum}
we know that $E_G \equiv m-n\; \text{mod}\;2$ and $O_G \equiv m\; \text{mod}\;2$. This implies that
$O_G = |T| \equiv m= |E|\; \text{mod}\;2$ and so
$|T|\equiv |E|\; \text{mod}\;2$, which is equivalent to $|T|+|E|$ is even.\\
On the other hand, fix some set of vertices $T$ such that $|T|\equiv |E|\; \text{mod}\;2$ and consider  a random orientation on $G$ obtained by directing each edge in $G$ independently of
the others and with probability $\frac{1}{2}$ in each direction. Let $E_G, O_G$ be as above.
We prove that there is a positive probability that the vertices of $T$ are the only ones
having odd degree. Since
$E_G \equiv m-n \equiv |T|-n\; \text{mod}\;2$ it follows that $n-|T|$ belongs to the range of $E_G$.
The result follows from Corollary \ref{fairdistr}, since
$\mathbb{P}[E_G =n-|T|]=\frac{1}{2^{n-1}}\binom{n}{n-|T|}>0$, and from the fact that any set, $T$, of $|T|\equiv |E|\; \text{mod}\;2$ vertices
can be such that all  vertices in $T$ have odd in-degree.
\end{proof}

We can also deduce a result on enumeration of oriented graphs.

\begin{lemma} Let $G=(V,E)$ be a graph on $n$ vertices and $m$ edges. Then the number of orientations
on the edges of $G$ for which there are exactly $t$ vertices of even in-degree equals $2^{m-n+1}\binom{n}{t}$.
\end{lemma}
\begin{proof} Let $\nu_t$ be the number of orientations of $G$ having exactly $t$ vertices of even in-degree.
Note that $t$ has to be such that $t\equiv m-n\; \text{mod}\; 2$.
From the set of all possible orientations of $G$, choose one uniformly
at random and let $A_t$ be the event that the orientation has $t$
vertices of even in-degree. Then
\[ \mathbb{P}[A]=\frac{\nu_t}{2^m} . \]
Now consider a random orientation on the edges of $G$ by directing each edge in $G$ independently of the others
and with probability $\frac{1}{2}$ in each direction. The result follows since,
by Corollary \ref{fairdistr},
the probability that there are $t$ vertices of even in-degree
equals $\frac{1}{2^{n-1}}\binom{n}{t}$.
\end{proof}

For similar results see \cite{Tazama}.
In a similar way, using Corollary \ref{randgraph}, one can obtain a result on enumeration of labeled graphs.
We leave the details to the reader.

\begin{lemma} The number of labeled graphs on $n$ vertices for which there are exactly
$t$ (where $t$ is even) vertices of odd degree equals $2^{m-n+1}\binom{n}{t}$.
\end{lemma}

Note that the case $t=0$ of the previous result appears as problem $16$ in $\S 5$ of \cite{Lovasz}.

\section{Open problems}\label{open}

There are many interesting questions concerning randomly oriented graphs.
So far we have studied the distribution of $E_G$, the number of vertices
with even in-degree. A natural generalization would be to consider the distribution of the number
of vertices whose in-degree equals $i \;\text{mod}\; k$. \\

Another interesting random variable is the number of vertices
with \emph{zero} in-degree. Suppose that we orient each edge of a connected graph $G$
independently and with probability $\frac{1}{2}$ for each direction. Let $Z_G$ be the number of
vertices with zero in-degree.
Thus, if
$Z_{G} = j$ then there are $j$ vertices in the
graph whose in-degree is zero.
Notice that the vertices of $G$ with zero in-degree form an independent set of vertices, i.e., no two of them are adjacent. Now we ask the following question. \\

For which graphs, $G$, is the distribution of $Z_G$ unimodal? \\

The distribution of $Z_G$ is related to the family of independent sets in $G$. If $Z_G=j$,
then $j$ vertices have in-degree zero and these $j$ vertices form an independent set.
That is, $Z_G=j$ gives rise to an independent set of vertices in $G$ of cardinality $j$.
So we might also ask the following. \\

For $j=0,1,\ldots,n$, denote by $\alpha_j(G)$ the
number of independent set of vertices of $G$ of cardinality $j$.
Is the sequence $\{\alpha_i(G)\}_{j=0}^{n}$ unimodal? \\

This problem is considered in \cite{AEMS}
where it is proven that the answer to the last question is no, for general graphs.
However, it remains an open question to determine whether the question is true in the case of trees.
In \cite{AEMS} one can find the following.

\begin{conj}[Alavi, Erd\H{o}s, Malde, Schwenk, $1987$] If $G$ is a tree, then the independent set sequence $\{\alpha_i(G)\}_{j=0}^{n}$ is unimodal.
\end{conj}

We believe that a similar result holds true for the distribution of $Z_G$, when $G$ is a tree.

\begin{conj} If $G$ is a tree, then the distribution of $Z_G$ is unimodal.
\end{conj}

\section*{Acknowledgement}
I am thankful to Robbert Fokkink and Tobias M\"{u}ller for many valuable discussions and comments.


\begin{thebibliography}{3}


\bibitem{AEMS} Yusef Alavi, Paul Erd\H{o}s, Paresh J. Malde and Allen J. Schwenk, \textit{The vertex independence sequence of a graph is not constrained}, Congr. Numer. \textbf{58} (1987), p 15--23.

\bibitem{AlmLinusson} Sven Erik Alm and Svante Linusson, \textit{A counter-intuitive correlation in a random tournament},
Combinatorics, Probability and Computing \textbf{20} (2011), no. 1, p. 1--9.

\bibitem{Bollobas} B\'ela Bollob\'as, \textit{Random Graphs}, Cambridge Studies in Advanced Mathematics, \textbf{73},
Cambridge University Press, 2001.

\bibitem{Meester} Erik Broman, Tim van de Brug, Wouter Kager and Ronald Meester, \textit{Stochastic
domination and weak convergence of conditioned Bernoulli random variables}, ALEA Lat. Am. J. Probab. Math. Stat. 9 (2012),
no. 2, p. 403--434.

\bibitem{Diarmir} Colin McDiarmir, \textit{General percolation and random graphs}, Adv. in Appl. Prob. \textbf{13} (1981), p. 40--60.

\bibitem{Frank} Andr\'{a}s Frank, Tibor Jord\'{a}n, Zolt\'{a}n Szigeti, \textit{An orientation theorem with parity conditions}, Discrete
Appl. Math. \textbf{115} (2001), p. 37--47.

\bibitem{Hamza} Kais Hamza, \textit{The smallest uniform upper bound on the distance between the mean
and the median of the Binomial and Poisson distributions}, Statistics \& Probability Letters \textbf{23} (1995), p. 21--25.


\bibitem{Hoeffding} Wassily Hoeffding, \textit{On the Distribution of the Number of Successes in Independent Trials},
An. Math. Statistics \textbf{27} (1956), p. 713-721.

\bibitem{Kaas} Rob Kaas and Jan M. Buhrman,  \textit{Mean, Median and Mode in Binomial Distributions},
Statistica Neerlandica \textbf{34}(1), (1980), p. 13--18.


\bibitem{Linusson} Svante Linusson, \textit{A note on correlations in randomly oriented graphs}, Preprint 2009,
arXiv:0905.2881.

\bibitem{Lovasz} L\'aszl\'o Lov\'asz, \textit{Combinatorial Problems and Excercises}, North-Holland, Amsterdam, 1979.

\bibitem{PelSch} Christos Pelekis and Moritz Schauer, \textit{Network coloring and colored coin games}, Chapter 4 in Search theory: A game-theoretic prespective, S. Alpern et al (eds.), 2013.

%\bibitem{Penrose} Mathew D. Penrose, \textit{Random Geometric Graphs}, volume 5 of Oxford Studies in Probability,
%Oxford University Press, 2003.


%\bibitem{Scott} Alex D. Scott, \textit{On induced subgraphs with all degrees odd}, Graphs and Combinatorics, 17, 2001,
%p. 539--553.

\bibitem{SS} Moshe Shaked and George J. Shanthikumar, \textit{Stochastic Orders and their Applications}.
Springer, New York, 2007.

\bibitem{Tazama} Shinsei Tazama, Teruhiro Shirakura and Saburo Tamura, \textit{Enumeration of digraphs with given number of vertices of odd out-degree and vertices of odd in-degree}, Discrete Mathematics \textbf{90} (1991), p. 63--74.


\bibitem{Wang} Yaohong H. Wang: \textit{On the number of successes in independent trials}, Statistics Sinica \textbf{3} (1993),
p. 295--312.

\end{thebibliography}
\end{document}